\documentclass[a4paper,11pt]{amsart}
\usepackage[english]{babel}
\usepackage{wrapfig}
\usepackage{graphicx}
\usepackage{subfigure}
\usepackage{amsmath}
\usepackage{amssymb}
\usepackage{amsfonts}
\usepackage{mathrsfs}
\usepackage{amsthm}
\usepackage{fancyhdr}
\usepackage{enumerate}
\usepackage{braket}
\usepackage{youngtab}
\usepackage{tikz}
\usepackage{listings}
\usetikzlibrary{snakes}
\usetikzlibrary{shapes.geometric}
\usepackage[all,cmtip,2cell]{xy}
\usepackage{color}   %May be necessary if you want to color links
\definecolor{darkgreen}{rgb}{0.0, 0.5, 0.0}

\usepackage{geometry}
\usepackage{color}

\theoremstyle{plain}                    %stile corsivo
\newtheorem{theo}{Theorem}[section]     %definizione ambiente teorema
\newtheorem{prop}[theo]{Proposition}    %definizione ambiente proposizione
\newtheorem{cor}[theo]{Corollary}       %definizione ambiente corollario
\theoremstyle{definition}              %stile roman
\newtheorem{defin}[theo]{Definition}    %definizione ambiente definizione
\theoremstyle{remark}                  %stile per osservazioni
\newtheorem{remark}[theo]{Remark}          %definizione ambiente osservazione

\newcommand{\z}{\mathbb{Z}}
\newcommand {\C}{\mathbb{C}}

\newcommand{\zz}{\mathbb{Z}/2\mathbb{Z}}
\newcommand{\px}{(p\in X)}

\newcommand{\yd}{(Y,D)}
\newcommand{\oyd}{(\bar{Y},\bar{D})}
\newcommand{\pyd}{(Y',D')}
\newcommand{\od}{\bar{D}}

\newcommand{\pl}{\mathbb{P}^1}

\newcommand{\aaut}{\overline{\textnormal{Aut}}}
\newcommand{\aut}{\textnormal{Aut}}

\title{Deformation types of Looijenga pairs of small length}
\author{Angelica Simonetti}
\date{}

\begin{document}

\begin{abstract}
Following the work already done by E. Looijenga and others, we provide an analysis of the deformation types of Looijenga (or anticanonical) pairs $\yd$, where the anticanonical divisor $D$ is made of $n$ irreducible components and $6\leq n \leq9$. In doing so we also give a description of their toric models.
\end{abstract}

\maketitle
\tableofcontents

\section*{Introduction}
Looijenga pairs have been studied over the years for their relation with cusp singularities, that is a specific type of isolated surface singularity, among other things. In 1981, Looijenga conjectured that a cusp singularity $\px$ was smoothable if and only if a certain Looijenga pair $\yd$ existed (\cite{L81}, see proposition (2.8) and the discussion in (2.9) or Corollary 2.2 in \cite{FM83} and the following conjecture for a more concise exposition of the subject). He also proved the necessity of this statement. To give a precise formulation of the conjecture we will need to introduce the concept of \emph{dual cusp}. Briefly, given a cusp $\px$ with minimal resolution $E\subset\tilde X$, the exceptional locus of the minimal resolution is either a rational node or a cycle of rational curves. To every cusp one can associate a dual cusp (see for example \cite{FM83} for more details). We will refer to the exceptional locus of this dual cusp as the \emph{dual cycle} to $\px$. Then we have:
\begin{theo}
The cusp $\px$ is smoothable if and only if the dual cycle $D$ sits as an anticanonical divisor on a smooth rational surface $Y$.
\end{theo}
The sufficiency of this theorem has been proved in 2015 by work of Gross, Hacking and Keel (\cite{GHK15a}), who obtain the proof of this statement as a consequence of a bigger mirror symmetry result, and then with a different approach by Engel (\cite{E15}). Three other insightful papers on cusps and Looijenga pairs are: an article by Friedman and Miranda on smmothability of cusps of small length (\cite{FM83}) that appeared before the two proofs mentioned above, a paper on cusps by Engel and Friedman (\cite{EF16}) and finally a quite comprehensive work on Looijenga pairs by Friedman alone (\cite{F15}).

A consequence of all this work is that we get a lower bound on the number of smoothing components of a cusp as stated in the proposition below, where $D$ is the cycle of curves dual to the cusp $\px$.

\begin{prop}[\cite{F15}, Theorem 3.16]
The number of smoothing components of the cusp singularity $\px$ with minimal resolution $E\subset\tilde X$ is greater than or equal to the number of deformation families of Looijenga pairs $\yd$.
\end{prop}

It has also been conjectured that there is in fact a bijection between the number of smoothing components of a cusp singularity and the number of deformation types of the correspondent Looijenga pairs (see the end of Example 4.5 in \cite{E15}).

A first analysis of the number of deformation types for Looijenga pairs of small length (here by \emph{length} of a Looijenga pair $\yd$, we mean the number of irreducible components of $D$) is contained in the paper by Looijenga already mentioned above (\cite{L81}, Chapter 1), where he shows that for each $n\leq5$ and each fixed $D$ there exists only one deformation type of pairs $\yd$. In this work we continue the analysis studying the deformation types of Looijenga pairs of length $n$, where $6 \leq n \leq 9$, obtaining the following result.

\begin{theo}[Theorem \ref{maint}]
    If $n=6,7$ or $n=8$ and $D$ has associated cycle of integers different from $(a,2,b,2,c,2,d,2)$ then there is one deformation type of negative definite Looijenga pairs $\yd$ of length $n$ with fixed $D$. If $n=8$ and $D$ is of type $(a,2,b,2,c,2,d,2)$ there are two deformation types, distinguished by $\pi_1(U)$, where $U=Y\setminus D$. Finally, if $n=9$, then there are at most three deformation types of negative definite Looijenga pairs $\yd$ of length $9$ with fixed $D$
\end{theo}
This theorem has been used in a paper by J. Li (\cite{Li22}) while proving that some special Looijenga pairs $\yd$ with split mixed Hodge structure and such that $D$ consists of six components are examples of Mori Dream Spaces.
The structure of this work is as follows. In the first section we present a brief introduction to Looijenga pairs in general with the main definitions and properties. In the second section we give a description of the toric models of Looijenga pairs of length $n$, where $6 \leq n \leq 9$. Finally section 3 contains the proof of the main result discussed above.

\subsection*{Acknowledgments}
The author wishes to thank her PhD advisor Paul Hacking for all the useful convesations, for his guidance and his patience. the author is also grateful to her mentor Jonny Evans for his support and advice.

\section{Loojenga pairs}

\begin{defin}
A \emph{Looijenga pair} or \emph{anticanonical pair} $\yd$ is a smooth projective surface $Y$ together with a connected singular nodal divisor $D\in |-K_Y|$ which is either an irreducible rational curve with a single node or a cycle of smooth rational curves, $D=\sum_{i=1}^n D_i$, where each $D_i$ meets $D_{i+1}$ transversally, with $i$ understood mod $n$.
\end{defin}

The integer $n$ is called the length of $D$, if the components of $D$ are indexed as above, we refer to $\yd$ as a \textit{labeled Looijenga pair} and to the sequence of self intersections $(-D_1^2,-D_2^2,\dots,-D_n^2)$ as the \textit{cycle of integers} associated to it. To fix the notation, we will always label the components of $D$ starting from the top-right one, for instance, for $n=6$ we have:
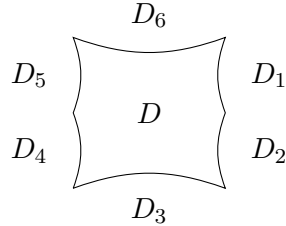
\begin{center}
\begin{tikzpicture}[out=20,in=160,relative]
\draw (-1,0) to (1,0) to (1,1) to (1,2) to (-1,2) to (-1,1) to  (-1,0);
\node [above] at (0,2) {$D_6$};
\node [right] at (1.2,1.5) {$D_1$};
\node [right] at (1.2,0.5) {$D_2$};
\node [below] at (0,0) {$D_3$};
\node [left] at (-1.2,0.5) {$D_4$};
\node [left] at (-1.2,1.5) {$D_5$};
\node at (0,1) {$D$};
\end{tikzpicture}
\end{center}
Note, as always, that all the pictures that will apear in this work are merely sketches: all components of $D$ should be understood as meeting transversally. An orientation of $D$ is an orientation of its dual graph, or equivalently the choice of a generator of $H_1(D,\z)\cong\z$. Observe that for $n\geq3$ an orientation determines a natural labeling of the components of $D$ up to cyclic permutation and viceversa a labeling induces an orientation on $D$.

\begin{defin}\label{autyd}
An \textit{isomorphism of labeled Looijenga pairs} $(Y,D)$ and $(Y',D')$ is an isomorphism $f:Y\rightarrow Y'$ such that $f(D_i)=D'_i$ for each $i=1,\dots,n$ which is compatible with the orientation of $D$ and $D'$. Let $\aut \yd$ be the group of automorphisms of a labeled Looijenga pair mapping each component of $D$ to itself and preserving the orientation of $D$.
\end{defin}

If the intersection matrix $(D_i\cdot D_j)$ is negative definite, we call $\yd$ a \textit{negative definite Looijenga pair} and say that $D$ is negative definite. A useful invariant of anticanonical pairs is their charge.
\begin{defin}[\cite{F15}, Definition 1.1]
The \textit{charge} $Q\yd$ of a Looijenga pair is defined as
$$Q\yd=12-D^2-n$$
\end{defin}
To give a glimpse of the cohomology theory of anticanonical pairs, let $\Lambda\yd\subset H^2(Y,\z)$ be the orthogonal complement of the lattice spanned by the classes of the $D_i$. Then $\Lambda\yd$ is free (\cite{F15}, Lemma 1.5) and, if $D$ is negative definite (which implies that the classes $D_i$ are independent in cohomology), its rank is equal to the charge minus two (\cite{F15}, Lemma 1.5). We also note that in the case $D$ is negative definite, then $Q\yd\geq 3$ (\cite{F15}, Corollary 1.3).\\
Always with the aim of fixing our notation let us give the following definitions:
\begin{defin}
Let $\yd$ be a Looijenga pair. A curve $C$ in $Y$ is an \textit{interior curve} if none of its irreducible components is contained in $D$. An \textit{internal (-2)-curve} instead is a smooth rational curve of self intersection -2 that is disjoint from $D$. We say that $\yd$ is \textit{generic} if it has no internal (-2)-curves.
\end{defin}
Define a \textit{simple toric blowup} to be the blowup of a Looijenga pair $\yd$ at a node of $D$ and an \textit{interior blowup} to be a blowup of $Y$ at a smooth point on $D$. For an interior blowup $\widetilde{Y}\rightarrow Y$, set $\widetilde{D}=\sum_i \widetilde{D}_i$, where $\widetilde{D}_i$ is the strict transform of $D_i$, while for a toric blowup define $\widetilde{D}=\sum_i \widetilde{D}_i+E$, where $\widetilde{D}_i$ is the strict transform of $D_i$ and $E$ is the exceptional divisor. Then in both cases $(\widetilde{Y},\widetilde{D})$ is still a Looijenga pair. Interior blowups increase the charge $Q\yd$ by one, while corner blowups do not change it (\cite{FM83}, Lemmas 3.3 and 3.4). Finally we observe that the charge of a Looijenga pair $\yd$ has a topological interpretation: let $U=Y\setminus D$, then $e(U)=Q(Y,D)$ where $e(U)$ is the Euler number of $U$ (\cite{F15}, Lemma 1.2).

\section{Toric models for Looijenga pairs of length $6\leq n\leq 9$}\label{toricmod}

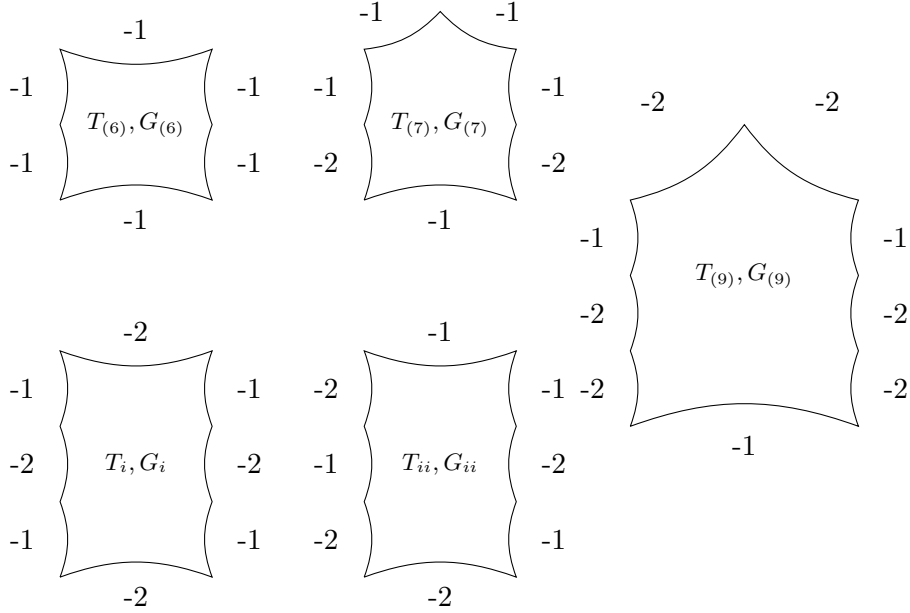
\begin{figure}[h]
    \begin{center}
    \begin{tikzpicture}[out=20,in=160,relative]
    \draw (-1,0) to (1,0) to (1,1) to (1,2) to  (-1,2) to (-1,1) to (-1,0);
    \node [above] at (0,2) {-1};
    \node [right] at (1.2,0.5) {-1};
    \node [right] at (1.2,1.5) {-1};
    \node [below] at (0,0) {-1};
    \node [left] at (-1.2,0.5) {-1};
    \node [left] at (-1.2,1.5) {-1};
    \node at (0,1) {\footnotesize{$T_{(6)},G_{(6)}$}};
    
    \draw (3,0) to (5,0) to (5,1) to (5,2) to (4,2.5) to  (3,2) to (3,1) to  (3,0);
    \node [right] at (4.6,2.5) {-1};
    \node [right] at (5.2,1.5) {-1};
    \node [right] at (5.2,0.5) {-2};
    \node [below] at (4,0) {-1};
    \node [left] at (2.8,0.5) {-2};
    \node [left] at (2.8,1.5) {-1};
    \node [left] at (3.4,2.5) {-1};
    \node at (4,1) {\footnotesize{$T_{(7)},G_{(7)}$}};
    
    \draw (-1,-5) to (1,-5) to (1,-4) to (1,-3) to (1,-2) to (-1,-2) to (-1,-3) to (-1,-4) to  (-1,-5);
    \node [above] at (0,-2) {-2};
    \node [right] at (1.2,-2.5) {-1};
    \node [right] at (1.2,-3.5) {-2};
    \node [right] at (1.2,-4.5) {-1};
    \node [below] at (0,-5) {-2};
    \node [left] at (-1.2,-4.5) {-1};
    \node [left] at (-1.2,-3.5) {-2};
    \node [left] at (-1.2,-2.5) {-1};
    \node at (0,-3.5) {\footnotesize{$T_{i},G_{i}$}};
    
    \draw (3,-5) to (5,-5) to (5,-4) to (5,-3) to (5,-2) to (3,-2) to (3,-3) to (3,-4) to (3,-5);
    \node [above] at (4,-2) {-1};
    \node [right] at (5.2,-2.5) {-1};
    \node [right] at (5.2,-3.5) {-2};
    \node [right] at (5.2,-4.5) {-1};
    \node [below] at (4,-5) {-2};
    \node [left] at (2.8,-4.5) {-2};
    \node [left] at (2.8,-3.5) {-1};
    \node [left] at (2.8,-2.5) {-2};
    \node at (4,-3.5) {\footnotesize{$T_{ii},G_{ii}$}};
    
    \draw (6.5,-3) to (9.5,-3) to (9.5,-2) to (9.5,-1) to (9.5,0) to (8,1)  to (6.5,0) to (6.5,-1) to (6.5,-2) to  (6.5,-3);
    \node [right] at (8.8,1.3) {-2};
    \node [right] at (9.7,-0.5) {-1};
    \node [right] at (9.7,-1.5) {-2};
    \node [right] at (9.7,-2.5) {-2};
    \node [below] at (8,-3) {-1};
    \node [left] at (6.3,-2.5) {-2};
    \node [left] at (6.3,-1.5) {-2};
    \node [left] at (6.3,-0.5) {-1};
    \node [left] at (7.1,1.3) {-2};
    \node at (8,-1) {\footnotesize{$T_{(9)},G_{(9)}$}};
    
    \end{tikzpicture}
    \end{center}
    \caption{Boundary cycles of the common toric pairs}\label{toricpairs}
    \end{figure}

Among Looijenga pairs there are some special ones which can be used to classify and analyze all the others, namely toric models and minimal pairs.
\begin{defin}
A Looijenga pair $\oyd$ is a \textit{toric pair} if $\bar{Y}$ is a smooth projective toric surface and $\od=\bar{Y}\setminus (\C^*)^2$ is the toric boundary. Now let $\pi:Y\rightarrow \bar{Y}$ be a sequence of interior blowups and let $D$ be the strict transform of $\od$: we call $\pi$ a \textit{toric model} for the Looijenga pair $\yd$.
\end{defin}

In other words, we say that $\yd$ admits a toric model if there exists a sequence of interior blow-downs $\pi:\yd\rightarrow\oyd$ where $\oyd$ is toric. We explicitly note that the charge of a toric pair is equal to zero.
\begin{remark}[\cite{GHK15b}]\label{toriciso}
The general theory of smooth projective toric varieties implies that the isomorphism type of a toric Looijenga pair is determined by its cycle of integers.
\end{remark}
Changing perspective, given any Looijenga pair we can always contract a sequence of $(-1)$-curves on it until we get to a pair $(Y',D')$, that we will call \textit{minimal}, where $Y'$ is a minimal rational surface. We have the following result by Miranda (\cite{M90}, Theorem 2.1). Note that from now until the end of the section we will assume that the divisor $D$ does not contain any $(-1)$-curves.
\begin{theo}
Let $\yd$ be a negative definite Looijenga pair with $D$ of length $n\geq 4$. Then $Y$ can be blown down to $\pl\times\pl$ so that $D$ is mapped to the standard square $D'=(\pl\times\{0,\infty\})\cup(\{0,\infty \}\times\pl)$.
\end{theo}

We observe that, in the hypothesis of the theorem above ($n\geq4)$, we can always arrange the sequence of blowups from $\pl\times\pl$ to $\yd$ so that we first perform all the toric blowups and then all the interior blowups, thus every negative definite anticanonical pair admits a map to a toric pair (with an exceptional cycle of the same length) which consists of a sequence of interior blowups, or, equivalently, every negative definite Looijenga pair admits a toric model. Let's focus on toric pairs for which $6\leq n\leq9$ (the cases for $n\leq 5$ have already been stuedied in \cite{L81}).

\begin{figure}[h]
    \begin{center}
    \begin{tikzpicture}[out=20,in=160,relative]
    \draw (-1,0) to (1,0) to (1,1) to (1,2) to (-1,2) to (-1,1) to  (-1,0);
    \node [above] at (0,2) {-1};
    \node [right] at (1.2,1.5) {0};
    \node [right] at (1.2,0.5) {-1};
    \node [below] at (0,0) {-1};
    \node [left] at (-1.2,0.5) {-2};
    \node [left] at (-1.2,1.5) {-1};
    \node at (0,1) {$\od$};
    \node at (3,1) {$\xymatrix{\ar[r]^{\varphi_{1,4}}&}$};
    \draw (5,0) to (7,0) to (7,1) to (7,2) to (5,2) to (5,1) to  (5,0);
    \node [above] at (6,2) {-1};
    \node [right] at (7.2,1.5) {-1};
    \node [right] at (7.2,0.5) {-1};
    \node [below] at (6,0) {-1};
    \node [left] at (4.8,0.5) {-1};
    \node [left] at (4.8,1.5) {-1};
    \node at (6,1) {$G_{(6)}$};
    \end{tikzpicture}
    \end{center}
    \caption{Case (b)-i}\label{b1}
    \end{figure}
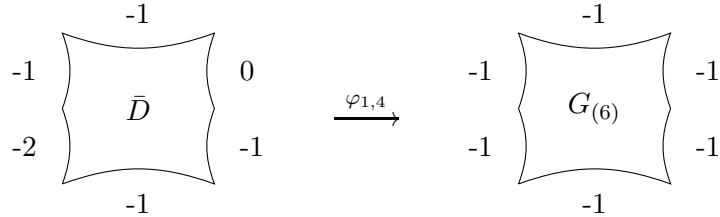

\begin{prop}\label{toricmodels}
Every negative definite Looijenga pair $\yd$ with cycle $D$ of length $n=6,7$ or $9$ can be blown down, for each $n$, to one common toric pair $(T_{(n)},G_{(n)})$ with $\textit{length}(D)=\textit{length}(G_{(n)})$. Looijenga pairs $\yd$ with length $n=8$ can always be blown down to one of the two toric pairs $(T_i,G_i)$ or $(T_{ii},G_{ii})$ whose toric boundaries are described in figure \ref{toricpairs}, along with the ones for $n=6,7,9$.
\end{prop}

% \begin{itemize}
% \item [i.] if $n=6$ every pair $\yd$ can be blown down to a single toric pair $(Y_,D_0)$ where $Y_0$ is a Del Pezzo surface of degree 6 and $D_0$ is the cycle of (-1)-curves contained in it.
% \item [ii.] if $n=7$ every pair $\yd$ can be blown down to a single toric pair $(Y_0,D_0)$ as well: it is given by the blow up of the toric model above in one point
% \item [iii.] if $n=8$ every pair can be blown down to either $(Y_1,D_1)$ or $(Y_2,D_2)$, the two toric pairs obtained by the one in (ii) via a single blow up
% \item [iv.] finally, if $n=9$ every pair $\yd$ can be blown down to a single toric pair $(Y_0,D_0)$: it is given by the blow up of the toric model $(Y_2,D_2)$ in one point
% \end{itemize}

\begin{proof}
We begin with the case where $n=6$. From the proof of Miranda's theorem we know that there always exists a ruling on $\yd$ such that two disjoint components of $D$ are sections of it (see the first \emph{claim} in the proof of Theorem 2.1 given in \cite{M90}, right after the first sentence). Thus, up to symmetry, there are two possible cases:
\begin{itemize}
\item[a] $D_6$ and $D_3$ are sections,
\item[b] $D_6$ and $D_2$ are sections.
\end{itemize}
Consider fibres of this ruling which do not contain any component of $D$: they are always chains of interior $(-2)$-curves with two $(-1)$-curves at the ends of the chain intersecting $D$. Indeed every negative definite Looijenga pair $\yd$ is obtained from a Looijenga pair $\oyd$ such that $\bar Y$ is a $\pl$-bundle, by a sequence of blowups. Therefore the fibres of the ruling on $\bar Y$ are smooth irreducible curves of self intersection 0. Since $\yd$ is obtained from $\oyd$ through a series of (either toric or interior) blowups, then the fibres of the ruling on $Y$ not containing any component of $D$ have to be chains of the type we described above. Fibres containing components of $D$ have a similar configuration: let $f=\cup F_i$ be such a fibre. Then some of the curves $F_i$ are components of $D$ and have no restrictions on their self intersections (other than the negative definiteness condition) while the others will be arranged in chains of $(-2)$-curves with a (-1)-curve at the end that intersects $D$.

%\begin{figure}[h]
%\begin{center}
%\begin{tikzpicture}[out=20,in=160,relative]
%\draw (-1,0) to (1,0) to (1,1) to (1,2) to (-1,2) to (-1,1) to  (-1,0);
%\node [above] at (0,2) {-1};
%\node [right] at (1.2,1.5) {0};
%\node [right] at (1.2,0.5) {-1};
%\node [below] at (0,0) {-1};
%\node [left] at (-1.2,0.5) {-2};
%\node [left] at (-1.2,1.5) {-1};
%\node at (0,1) {$\od$};
%\node at (3,1) {$\xymatrix{\ar[r]^{\varphi_{3,6}}&}$};
%\draw (5,0) to (7,0) to (7,1) to (7,2) to (5,2) to (5,1) to  (5,0);
%\node [above] at (6,2) {-1};
%\node [right] at (7.2,1.5) {-1};
%\node [right] at (7.2,0.5) {-1};
%\node [below] at (6,0) {-1};
%\node [left] at (4.8,0.5) {-1};
%\node [left] at (4.8,1.5) {-1};
%\node at (6,1) {$D'$};
%\end{tikzpicture}
%\end{center}
%\caption{Case (a)}\label{a}
%\end{figure}
\begin{figure}[h]
    \begin{center}
    \begin{tikzpicture}[out=20,in=160,relative, scale=1.00]
    \draw (-1,0) to (1,0) to (1,1) to (1,2) to (-1,2) to (-1,1) to  (-1,0);
    \node [above] at (0,2) {-1};
    \node [right] at (1.2,1.5) {0};
    \node [right] at (1.2,0.5) {0};
    \node [below] at (0,0) {-2};
    \node [left] at (-1.2,0.5) {-1};
    \node [left] at (-1.2,1.5) {-2};
    \node at (0,1) {$\od$};
    \node at (2.8,1) {$\xymatrix{\ar[r]^{\varphi_{1,3}}&}$};
    \draw (4.6,0) to (6.6,0) to (6.6,1) to (6.6,2) to (4.6,2) to (4.6,1) to  (4.6,0);
    \node [above] at (5.6,2) {-1};
    \node [right] at (6.8,1.5) {0};
    \node [right] at (6.8,0.5) {-1};
    \node [below] at (5.6,0) {-1};
    \node [left] at (4.4,0.5) {-2};
    \node [left] at (4.4,1.5) {-1};
    \node at (5.6,1) {$D'$};
    \node at (8.4,1) {$\xymatrix{\ar[r]^{\varphi_{2,5}}&}$};
    \draw (10.2,0) to (12.2,0) to (12.2,1) to (12.2,2) to (10.2,2) to (10.2,1) to  (10.2,0);
    \node [above] at (11.2,2) {-1};
    \node [right] at (12.4,1.5) {-1};
    \node [right] at (12.4,0.5) {-1};
    \node [below] at (11.2,0) {-1};
    \node [left] at (10,0.5) {-1};
    \node [left] at (10,1.5) {-1};
    \node at (11.2,1) {$D''$};
    \end{tikzpicture}
    \end{center}
    \caption{Case (b)-ii}\label{b2}
\end{figure}
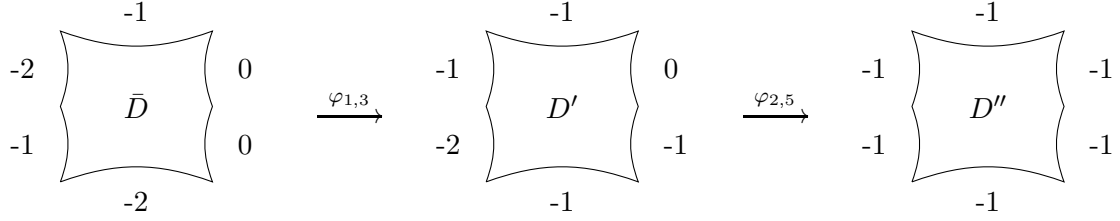

Now we can contract all singular fibres in the ruling which do not contain any components of $D$ (always on the most negative section) until they are irreducible fibres and blow down singular fibres containing components of $D$ to a chain contained in $D$. We get a map to a new pair $\pi:\yd\rightarrow\oyd$ where $\bar{Y}$ is a toric surface and $\od=\pi(D)$ its toric boundary. Let $\od_i$, for $i=1,\dots,6$, be the irreducible components of $\od$.\par
Suppose we are in case (a). Then $\od_1+\od_2$ and $\od_4+\od_5$ are fibers of the ruling on $\oyd$ induced by the one on $\yd$, therefore they must be pairs of (-1)-curves. Moreover, the fact that $\oyd$ is a toric pair implies that 
\begin{equation}\label{formula}
-\sum \od_i ^2=3\cdot n-12
\end{equation}
hence, in our case, $-\sum \od_i ^2=6$, so that $-\od_6^2-\od_3^2=2$. Furthermore, we know that $-D_3^2\geq 2$ and $-D_6^2\geq 2$, and by the algorithm we used we must have $|\od_3^2-\od_6^2|\leq 1$, therefore $\od_6^2=-1=\od_3^2$. In this case remark \ref{toriciso} implies that $\oyd$ is isomorphic to the pair $(T_{(6)},G_{(6)})$, where $T_{(6)}$ is a Del Pezzo surface of degree 6 and $G_{(6)}$ is the cycle of (-1)-curves contained in it. 

%Otherwise, without loss of generality, suppose we have $\od_6^2=-2$ and $\od_3^2=0$ (the case $\od_6^2=0$ and $\od_3^2=-2$ is treated analogously).

Next, suppose we are in case (b). Then $\od_1$ has to be a simple fibre with self intersection equal to zero. Similarly, $\od_3+\od_4+\od_5$ is a singular fibre and there are only two possible arrangements of self intersections for this triple of curves:
\begin{itemize}
\item [i.] $\od_3^2=-1$, $\od_4^2=-2$ and $\od_5^2=-1$
\item [ii.] $\od_3^2=-2$, $\od_4^2=-1$ and $\od_5^2=-2$
\end{itemize}
Let us start with (i). Using (\ref{formula}) we get that $-\od_6^2-\od_3^2=2$, thus, using the argument given for case (a), we can assume that $-\od_6^2=-1=-\od_3^2$ and $\oyd$ has associated cycle of integers $(0,1,1,2,1,1)$.

Now we observe (figure \ref{b1}) that there exists a ruling on $\oyd$ with sections given by the curves $\od_1$ and $\od_4$ with self intersections respectively 0 and -2.
Since the cycle $D$ we started with is negative definite (with no (-1)-curves on it), then it follows that $D_1^2$ and $D_4^2$ are both less or equal to -2. This means that when blowing down on singular fibres meeting $D_1$ and $D_4$ we contracted (at least) two (-1)-curves intersecting $D_1$ and none intersecting $D_4$, and we can always change that and contract one curve on $D_4$ and one on $D_1$. In other words, there exist an \textit{elementary transformation} $\varphi$
 $$\xymatrix{ & (\widehat{Y},\widehat{D})\ar[dl]_{\textnormal{Bl}_p}\ar[dr]^{\textnormal{Bl}_q} & \\ \oyd\ar@{-->}[rr]^{\varphi} & & (Y',D')}$$
given by the composition of the blowup of a smooth point $p$ on $\od_1$ and the blowdown of a (-1)-curve intersecting $\widehat{D}_4$, the strict transform of $\od_4$, to a smooth point $q$ on $D'_4$, and a map $\psi:\yd\rightarrow (\widehat{Y},\widehat{D})$ such that $\textnormal{Bl}_p\circ\psi=\pi$. This gives us a new toric model for the anticanonical pair we started with, $\yd\rightarrow (Y',D')$ and the diagram:
$$\xymatrix{\yd\ar[dd]_{\pi}\ar[dr]^{\psi}& & \\ & (\widehat{Y},\widehat{D})\ar[dl]_{\textnormal{Bl}_p}\ar[dr]^{\textnormal{Bl}_q} & \\ \oyd\ar@{-->}[rr]^{\varphi} & & (Y',D')}$$
where $D'$ is a cycle of six (-1)-curves (figure \ref{b1}). As a consequence, the toric pair $(Y',D')$ is again isomorphic to the anticanonical toric pair $(T_{(6)},G_{(6)})$.
\begin{remark}
From now on let us denote by $\varphi_{i,j}$ the elementary transformation from $\oyd$ to $(Y',D')$ that is given by the composition of the blowup of a smooth point $p$ on $\od_i$ and the blowdown of a (-1)-curve intersecting $\widehat{D}_j$, the strict transform of $\od_j$, to a smooth point $q$ on $D'_j$ (when it exists, i.e., when there is a ruling of $\oyd$ s.t. $\od_i$ and $\od_j$ are sections). For the birational map $\varphi$ described above, we would have $\varphi=\varphi_{1,4}$.
\end{remark}

% Again, since $D_1=0> -2$, there exists an elementary transformation $\varphi_{1,4}$ such that $\varphi_{1,4}(\od)$ is equal to a cycle of six (-1)-curves and we get the toric model $\yd\rightarrow (T_{(6)},G_{(6)})$.\\

Finally, in case (ii), we have $-\od_0^2-\od_3^2=1$ thus $\od$ is associated to the cycle of integers $(0,0,2,1,2,1)$ (see figure \ref{b2}). First we use the existence of the ruling with sections given by $\od_1,\od_3$ to get to the pair $(Y',D')$ with $D'$ given by the cycle of integers $(1,0,1,1,2,1)$ through the elementary transformation $\varphi_{1,3}$. Then we use the ruling with sections given by $\od_2,\od_5$ to get to the pair $(Y'',D'')$ with $D''$ given by the cycle of integers $(1,1,1,1,1,1)$.  We have obtained again a map made of of interior blow-downs from $\yd$ to $(T_{(6)},G_{(6)})$, thus proving the proposition for $n=6$.\par

For $n=7,8,9$ we proceed similarly, proving the required statement case by case. First we distinguish different cases, based on the list of possible arrangement of the two sections contained in $D$, up to symmetry. Then for each case we suppose to contract the singular fibres of the given ruling for $\yd$ following the procedure discussed above until we end up with a toric pair $\oyd$ of the same length. We make a list of all such toric pairs, depending on the possible arrangements of fibers and using (\ref{formula}) to determine the cycles of self intersections. Finally, for all  toric pairs that are different from the ones described in figure \ref{toricpairs} we use a sequence of elementary transformations to obtain the toric model we are looking for. For $n=7$ the common toric pair $(T_{(7)},G_{(7)})$ is obtained from $(T_{(6)},G_{(6)})$ through a toric blowup. There are two possible choices for the positions of the two sections in $D$ and each of these admits different toric pairs, depending on the self intersection of the fibres: these toric pairs and the corresponding elementary transformations are described in table \ref{table7}. For $n=8$ we fix two toric pairs which are, again, obtained from $(T_{(7)},G_{(7)})$ via one appropriate toric blowup: $(T_{i},G_{i})$ has cycle of integers $(1,2,1,2,1,2,1,2)$ and $(T_{ii},G_{ii})$ is associated to $(1,2,1,2,2,1,2,1)$. In this case the possible configurations of the sections in $D$ are three: for each of them table \ref{table8} lists all the toric pairs we could get and the relative elementary transformations used to obtain either $(T_{i},G_{i})$ or $(T_{ii},G_{ii})$. 
Finally, for $n=9$, $(T_{(9)},G_{(9)})$ is constructed from $(T_{ii},G_{ii})$ through a toric blowup in such a way that the correspondent cycle is given by $(2,2,1,2,2,1,2,2,1)$. Here there are again three possible configurations of sections contained in $D$. The list of toric pairs and relative elementary transformations can be found in table \ref{table9}, in the appendix.
\end{proof}

\begin{remark}\label{abrem}
Notice that the pairs $a,b$ in table \ref{table8} are associated to the same ruling and share the same arrangement of fibers: $\od_1,\od_2,\od_3$ with self intersections $-1,-2,-1$ and $\od_5,\od_6,\od_7$ again with self intersections $-1,-2,-1$. The only differences between them are the self intersections of the two sections. This is to take into account the negative definite Looijenga pairs $\yd$ obtained from $a$ if none of the interior blowups are performed on either $\od_4$ or $\od_8$. Indeed in this case there will be no reducible fibres intersecting $D_4$ and $D_8$, who will still have self intersections equal to $-2$. As a consequence such pairs can only be contracted to the toric pair $a$ and not to the toric pair $b$ (in other words there does not exist an elementary transformation from pair $a$ to pair $b$). The same reasoning applies to pairs $c,d$. We observe explicitly that $a$ and $d$ are the only two pairs that admit a sequence of elementary transformations taking them to $(T_{i},G_{i})$ but not one to $(T_{ii},G_{ii})$, while for all the other toric pairs listed in table \ref{table8} there exists a sequence of maps $\varphi_{i,j}$ bringing them to $(T_{i},G_{i})$, $(T_{ii},G_{ii})$ or even both of them.
\end{remark}

\begin{remark}\label{ellipticpairs} 
A direct consequence of the result above is that every non negative Looijenga pair with cycle $D$ of length $n=6,7,9$ can be contracted to a pair $(Y_{(n)}, D_{(n)})$ that is obtained from the toric pair $(T_{(n)},G_{(n)})$ of the appropriate length performing an interior blowup on every component of self intersection -1. Similarly, Looijenga pairs with cycle of length $n=8$ can always be contracted to either one of the anticanonical pairs $(Y_{i},D_{i})$ or $(Y_{ii},D_{ii})$, which again are obtained respectively from $(T_{i},G_{i})$ and $(T_{ii},G_{ii})$ through four interior blowups on the (-1)-curves. For ease of notation let us assume, without loss of generality, that these blowups are performed on components $2,3,4,5,7$ for $n=7$, $1,3,5,7$ for $D_i$, $1,3,6,8$ for $D_{ii}$ and $1,4,7$ for $n=9$. Note that, since all the irreducible components of the toric boundaries described in figure \ref{toricpairs} have self intersections less than or equal to 2 in absolute value, then such a pair $(Y_{(n)},D_{(n)})$ can always be constructed: from now on we will refer to these pairs as \textit{elliptic pairs}, since they are deformation euivalent to elliptic surfaces and their anticanonical cycles are singular fibres of the corresponding elliptic fibration.

Let us focus on the two elliptic pairs with cycles of length 8. Clearly the cycles of integers for these two pairs coincide, since both $D_{i}$ and $D_{ii}$ are made of eight curves of self intersection $(-2)$, thus we may ask ourselves if these anticanonical pairs are isomorphic or not. Observe that they can be distinguished by the fundamental group of the complement of their anticanonical divisors, $U_i:=Y_i\setminus D_i$ and $U_{ii}:=Y_{ii}\setminus D_{ii}$. More precisely, $\pi_1(U_i)=\zz$ while $\pi_1(U_{ii})=0$. This can be seen recalling that if $\yd$ is obtained from a toric surface through a sequence of interior blowups, then $\pi_1(Y\setminus D)=N/\langle v_1,\dots,v_p\rangle$, where $N$ is the lattice containing the fan of the toric variety and $v_1,\dots,v_p\in N$ are the primitive vectors corresponding to the curves of the toric boundary where the blowups are performed. Now, for $Y_{ii}$, the set of vectors $\{w_1,\dots,w_4\}$ corresponding to the four (-1)-curves contains a basis for the lattice $N$, thus the fundamental group of $U_{ii}$ is trivial, while the set of vectors $\{w_1,\dots,w_4\}$ in the fan associated to $T_i$ corresponding to the four (-1)-curves share the linear relations $v_1=-v_3,v_2=-v_4,v_1+v_2=2e_1,v_1+v_4=2e_2$ where $e_1,e_2$ are generators for the lattice $N$, thus $\pi_1(U_i)=\zz$.
\end{remark}

This remark, together with the proof of proposition \ref{toricmodels} allows us to refine the result of that proposition for Looijenga pairs of length eight.

\begin{cor}\label{toricmodels8}
Let $\yd$ be a negative definite Looijenga pair with anticanonical cycle of length $n=8$. Then either $\pi_1(Y\setminus D)=\zz$ or $\pi_1(Y\setminus D)=0$. If $\pi_1(Y\setminus D)=\zz$ the pair $\yd$ can only be contracted to the toric pair $(T_i,G_i)$, otherwise if $\pi_1(Y\setminus D)$ is trivial it can be contracted to both $(T_{ii},G_{ii})$ and $(T_i,G_i)$.
\end{cor}

\begin{proof}
Let $\yd$ be such that $\pi_1(Y\setminus D)=\zz$. Then $\yd$ must be obtained from the toric model $(T_i,G_i)$, since, as we saw in remark \ref{ellipticpairs}, this is the only way to keep the fundamental group non trivial. More precisely, the Looijenga pair $\yd$ is obtained from $(Y_i,D_i)$ blowing up points only on the odd components of $D_i$. Observe now that the fact that all the interior blowups happen on the odd components of $D_i$ implies that the fundamental group of $Y\setminus D$ is equal to the fundamental group of $(Y_i\setminus D_i)=\zz$. Thus, $\pi_1(Y\setminus D)=\zz$ if and only if $\yd$ is of the kind described above, in which case it only admits a birational map to $(T_i,G_i)$ On the other hand, if $\yd$ is a negative definite Looijenga pair of length eight with trivial fundamental group, then it can be contracted to either of the toric models listed in proposition \ref{toricmodels} (see table \ref{table8}).
\end{proof}

\section{Deformation types of Loojienga pairs for $6 \leq n\leq 9$}
Proposition \ref{toricmodels} allows us to give a complete description of the number of deformation classes of Looijenga pairs $\yd$ with fixed $D$ of length $6\leq n \leq 9$. In order to do that we need to recall briefly the description of the Mordell-Weil group of the surfaces $(Y_{(n)},D_{(n)})$ and $(Y_{i},D_{i}),(Y_{ii},D_{ii})$ constructed in remark \ref{ellipticpairs}, and more generally on their automoprhism groups. Given a rational elliptic surface $Y$ with section, the Mordell-Weil group can be thought of as a subgroup of the automorphism group of the surface itself and it acts transitively on sections of $Y\rightarrow \pl$: each element of the Mordell-Weil group gives an automorphism of $Y$ that acts as a translation on (smooth) fibres. Thus if $(Y,D)$ is an elliptic Looijenga pair, this group can also be identified with a subgroup of the generalized automorphism group $\aaut(Y,D)$ containing all automorphisms of $(Y,D)$ fixing $D$ set-wise but not component wise. Note that the automorphism group of $(Y,D)$ as defined in \ref{autyd} is also contained in $\aaut(Y,D)$ as a subgroup, more precisely it can be thought of as the kernel of the map $\gamma:\aaut(Y,D)\rightarrow \mathcal D_n$, with $n$ equal to the length of $D$, that sends each automorphism $\phi$ in $\aaut(Y,D)$ to the element of the dihedral group which corresponds to the action of $\phi$ to the dual graph for $D$. Hence we get the sequence of maps:
\begin{equation}\label{autseq}
0 \rightarrow \aut(Y,D)\rightarrow \aaut(Y,D)\xrightarrow{\gamma} \mathcal D_n\end{equation}
Now for the elliptic Looijenga pairs of remark \ref{ellipticpairs}, we can prove the following result.

\begin{prop}\label{autelliptic}
Let $(Y_{(n)},D_{(n)})$ be the elliptic Looijenga pair of length $n=6,7$. Then the automophism group $\aaut(Y_{(n)},D_{(n)})$ projects onto the dihedral group $\mathcal D_n$ of order $2n$. Similarly the automorphism group $\aaut(Y_{ii}, D_{ii})$ admits a sujective map onto $\mathcal D_8$. Instead, the automorphism group $\aaut(Y_{i},D_{i})$ admits a surjective map to $\mathcal D_4$, but not onto $\mathcal D_8$. The automorphism group $\aaut(Y_{(9)},D_{(9)})$ admits a surjective map to $\mathcal D_3$, but not onto $\mathcal D_9$.
\end{prop}
\begin{proof}
The proof will proceed as follows: for each $n$ listed above, first we will show that the Mordel Weil group admits a surjective map onto $\z/n\z$, then we will construct explicitly an involution of the elliptic anticanonical pair thus showing that the map $\gamma$ in sequence \ref{autseq} is surjective.\\
For $n=6$, then there is a (-1)-curve intersecting every component of $D_{(6)}$. Since each of these (-1)-curves is a section for the elliptic surface, and $\textnormal{MW}(Y_{(6)})$ acts transitively on sections, then there must esists an automorphism of $Y_{(6)}$ mapping $D_{(6),i}\mapsto D_{(6),i+1}$ and thus acting as a rotation of order 6 on the dual graph. We therefore get a surjective map $\textnormal{MW}(Y_{(6)})\rightarrow \z/6\z$. As for the involution, let us describe the toric pair of length 6 through its fan. Let $N$ be the lattice isomorphic to $\z^2$ and let $e_1,e_2$ be the vectors corresponding to $(1,0),(0,1)$. Then the cones of the fan for $(T_{(6)},G_{(6)})$ are generated by the vectors $\{e_1,-e_2\}$,$\{-e_2,-e_1-e_2\}$,$\{-e_1-e_2,-e_1\}$,$\{-e_1,e_2\}$,$\{e_2,e_1+e_2\}$ and $\{e_1+e_2,e_1\}$. The lattice isomorphism mapping $e_1\mapsto e_2$ induces an involution of the toric Looijenga pair which lifts to an involution of $(Y_{(6)}, D_{(6)})$ if the 6 interior blowups are made at the appropriate points. Since this involution acts on the dual graph of $D_{(6)}$ as a reflection our claim is proven.

Similarly, if $n=7$, since there is a curve intersecting the second and third components of $D_{(7)}$ (see remark \ref{ellipticpairs}) then again there exists an automorphism of the elliptic pair of length 7 mapping $D_{(7),i}\mapsto D_{(7),i+1}$ and thus acting as a rotation of order 7 on the dual graph. The fan for the toric pair of length 7 is obtained from the one above adding the ray generated by $-e_1+e_2$ and the required involution is obtained as a lift of the one induced on $(T_{7},G_{(7)})$ by the linear map sending $e_2$ to $-e_1$ (and fixing the new ray). The case of $(Y_{ii}, D_{ii})$ is treated similarly noticing that components $1,8$ are adjacent and both intersect a (-1)-curve given by an interior blowup. Moreover the linear map used for the previous case still gives us an involution of $(Y_{ii},D_{ii})$ following the same procedure as above once we observe that the toric pair $(T_{ii},G_{ii})$ can be constructed adding the rays $e_1+2e_2,-2e_1-e_2$ to the fan we gave for the case $n=6$.

Now consider $(Y_{i},D_{i})$. Given the arrangement of the interior blowups, the rotation of the dual graph with greatest order is the one induced by the automorphism mapping $D_{i,j}$ to $D_{i,j+2}$ which has order 4. Therefore we get a surjective map $\textnormal{MW}(Y_{i})\rightarrow \mathcal \z/4\z$, but not one onto $\mathcal \z/8\z$ because $\textnormal{MW}(Y_{i})\cong \z/4\z$ (see \cite{M89}, page 82). An involution can be constructed as usual, once we notice that the fan for the toric pair of length eight $(T_{i},G_{i})$ can be obtained from the one of the toric pair of length seven adding the ray $e_1-e_2$ and using the same linear map.

Finally, consider the pair  $(Y_{(9)},D_{(9)})$. Again, because of the arrangement of the interior blowups, the rotation of the dual graph with greatest order is the one induced by the automorphism mapping $D_{i,j}$ to $D_{i,j+3}$ which has order 4. Therefore we get a surjective map $\textnormal{MW}(Y_{(9)})\rightarrow \mathcal \z/3\z$, but not one onto $\mathcal \z/9\z$, again because $\textnormal{MW}(Y_{(9)})\cong \z/3\z$ (see \cite{M89}, page 82). Moreover, an involution can be constructed as usual, once we notice that the fan for the toric pair of length nine $(T_{(9)},G_{(9)})$ can be obtained from the one of the toric pair $(T_{ii},G_{ii})$ adding the ray $e_1-e_2$ and using the same linear map.
\end{proof}

\begin{theo}\label{maint}
If $n=6,7$ or $n=8$ and $D$ has associated cycle of integers different from $(a,2,b,2,c,2,d,2)$ then there is one deformation type of negative definite Looijenga pairs $\yd$ of length $n$ with fixed $D$. If $n=8$ and $D$ is of type $(a,2,b,2,c,2,d,2)$ there are two deformation types, distinguished by $\pi_1(U)$, where $U=Y\setminus D$. Finally, if $n=9$, then there are at most three deformation types of negative definite Looijenga pairs $\yd$ of length $9$ with fixed $D$.
\end{theo}
\begin{proof}
Let $\yd$, $\pyd$ be negative definite Looijenga pairs of length $n=6,7$. Thanks to remark \ref{ellipticpairs}, we know that they can always be obtained from the elliptic pair of the same length pair through a sequence of $m$ interior blowups, $\pi: \yd\rightarrow (Y_{(n)},D_{(n)})$ and $\pi':\pyd\rightarrow (Y_{(n)},D_{(n)})$. Let us fix labelings on $D,D',D_{(n)}$ so that $\pi(D_i)=D_{(n),i}$ and similarly $\pi'(D'_i)=D_{(n),i}$. Finally let us assume that they share the same cycle of integers $(a_1,\dots,a_n)$.  Then there must be an element $\sigma$ of the dihedral group of order $2n$ such that $(-D^2_{1} -2,\dots,-D^2_{n} -2)=(-D'^2_{\sigma(1)} -2,\dots,-D'^2_{\sigma(n)} -2)$. If $\sigma$ is the identity map, then all the blowups happen on the same components of $D_{(n)}$ for both $\yd$ and $\pyd$, therefore they are deformation equivalent. Otherwise, proposition \ref{autelliptic} guarantees that there exists an automorphism $\phi$ of the elliptic Looijenga pair of length $n$ which maps $D_{(n),i}\mapsto D_{(n),\sigma(i)}$, therefore $\yd$ and $\pyd$ are obtained from isomorphic pairs by performing the blowups on components of the anticanonical cycle which are identified via the isomorphism $\phi$ and again they are deformation equivalent.

Now let us assume that $\yd$, $\pyd$ are negative definite Looijenga pairs of length eight with fixed cycle of integers $(a_1,\dots,a_8)$. If $\pi_1(Y\setminus D)=\pi_1(Y'\setminus D')=0$, then both pairs can be contracted to $(Y_{ii},D_{ii})$ and thanks to proposition \ref{autelliptic} the same argument used above will show that the two Looijenga pairs are deformation equivalent. If $\pi_1(Y\setminus D)=\pi_1(Y'\setminus D')=\zz$, then it follows from corollary \ref{toricmodels8} that $\yd$ and $\pyd$ can only be contracted to the elliptic pair $(Y_i,D_i)$ and the blowups happen only on the odd components of $D_i$. Since $\aut(Y_i,D_i)$ admits a surjective map to $\mathcal D_4$ (proposition \ref{autelliptic}), then again the argument used in the previous cases allows us to show that $\yd$ and $\pyd$ are deformation equivalent. Finally let us suppose $\pi_1(Y\setminus D)=\zz$ and $\pi_1(Y'\setminus D')=0$ and suppose that they are deformation equivalent. Then there would exists a diffeomorphism mapping $\yd$ to $\pyd$ and sending $D$ to $D'$: this would imply in particular that the fundamental groups of $Y\setminus D$ and $Y'\setminus D'$ are isomorphic, thus giving a contradiction. Thus they cannot be deformation equivalent and the statement is proved.

If $n=9$, then given we only have a surjective map from the elliptic pair $\aut(Y_{(9)},D_{(9)})$ to $\mathcal D_3$, there are at most three choices for the arrangement of the interior blowups which do not need to give an automorphism of the elliptic pair. Therefore there are at most three deformation types.
\end{proof}

\medskip
\bibliography{C:/Users/agsim/Documents/Bibliography/biblio}
\bibliographystyle{alpha}

\appendix
\section{Tables, proof of theorem \ref{toricmodels}}
We include here the tables with all elementary transformations used in the proof of theorem \ref{toricmodels}.
\begin{table}[h]
    \begin{center}
    \begin{tabular}{ |c|c| } 
     \hline
     \multicolumn{2}{|c|}{Sections: $D_3$ and $D_7$} \\
     \hline
     Toric pair, cycle of integers & Elementary transformations \\ \hline
      (1,1,1,1,2,1,2) & Not needed  \\ 
     \hline
     \multicolumn{2}{|c|}{Sections: $D_2$ and $D_7$} \\
     \hline
     Toric pair, cycle of integers & Elementary transformations \\ \hline
      (0,1,1,2,2,1,2) & $\varphi_{4,1}$ \\ 
      (0,1,1,3,1,2,1) & $\varphi_{4,1}$ \\ 
      (0,0,2,2,1,3,1) & $\varphi_{3,1},\varphi_{6,2}$ \\ 
     \hline
    \end{tabular}
    \end{center}\caption{Elementary transformations for $n=7$}\label{table7}
    \end{table}

\begin{table}[h]
    \begin{center}
    \begin{tabular}{ |c|c|c| } 
     \hline
     \multicolumn{3}{|c|}{Sections: $D_4$ and $D_8$} \\
     \hline
     Toric pair, cycle of integers & Elementary transformations & Toric model \\ \hline
     $a.$ (1,2,1,2,1,2,1,2) & Not needed & $(T_{i},G_{i})$ \\ 
     $b.$ (1,2,1,1,1,2,1,3) & $\varphi_{6,3},\varphi_{8,5}$ & $(T_{ii},G_{ii})$ \\ 
     (2,1,2,1,1,2,1,2) & Not needed & $(T_{ii},G_{ii})$ \\
     \hline
     \multicolumn{3}{|c|}{Sections: $D_3$ and $D_8$} \\
     \hline
     Toric pair, cycle of integers & Elementary transformations & Toric model \\ \hline
     (1,1,2,1,2,2,1,2) & Not needed & $(T_{ii},G_{ii})$ \\
     (1,1,1,1,2,2,1,3) & $\varphi_{8,3}$ & $(T_{ii},G_{ii})$ \\
     \hline
     \multicolumn{3}{|c|}{Sections: $D_2$ and $D_8$} \\
     \hline
     Toric pair, cycle of integers & Elementary transformations & Toric model \\ \hline
     $c.$ (0,1,1,2,2,2,1,3) & $\varphi_{4,1},\varphi_{8,3}$ & $(T_{ii},G_{ii})$ \\ 
     $d.$ (0,2,1,2,2,2,1,2) & $\varphi_{5,1}$ & $(T_{i},G_{i})$ \\
     (0,1,1,2,3,1,2,2) & $\varphi_{8,2},\varphi_{5,1}$ & $(T_{ii},G_{ii})$ \\ 
     (0,1,1,3,1,3,1,2) & $\varphi_{4,1},\varphi_{4,1},\varphi_{6,3}$ & $(T_{ii},G_{ii})$ \\ 
     (0,1,1,3,2,1,3,1) & $\varphi_{4,1},\varphi_{7,2}$ & $(T_{ii},G_{ii})$ \\
     (0,1,1,4,1,2,2,1) & $\varphi_{4,1},\varphi_{4,1}$ & $(T_{ii},G_{ii})$ \\ 
     (0,1,2,1,4,1,2,1) & $\varphi_{5,1},\varphi_{3,8},\varphi_{5,2}$ & $(T_{ii},G_{ii})$ \\ 
     (0,1,2,2,2,1,4,0) & $\varphi_{7,1},\varphi_{3,8},\varphi_{7,2}$ & $(T_{ii},G_{ii})$ \\ 
     (0,1,2,3,1,2,3,0) & $\varphi_{7,1},\varphi_{4,1}$ & $(T_{ii},G_{ii})$ \\ 
     (0,1,3,1,3,1,3,0) & $\varphi_{7,1},\varphi_{7,1},\varphi_{4,8}$ & $(T_{ii},G_{ii})$ \\ 
     \hline
    \end{tabular}
    \end{center}\caption{Elementary transformations for $n=8$}\label{table8}
    \end{table}
    
    \begin{table}[h]
    \begin{center}
    \begin{tabular}{ |c|c| } 
     \hline
     \multicolumn{2}{|c|}{Sections: $D_1$ and $D5$} \\
     \hline
     Toric pair, cycle of integers & Elementary transformations \\ \hline
     (i) (2,1,2,1,3,1,2,2,1) & $\varphi_{7,2},\varphi_{1,6}, \varphi_{9,5}$\\ 
     (3,1,2,1,2,1,2,2,1) & $\varphi_{8,4}$, back to (i) \\ 
     (2,1,2,1,2,1,3,1,2) & $\varphi_{9,4}$, back to (i) \\
     (ii) (1,1,2,1,2,2,2,1,3) & $\varphi_{1,6}, \varphi_{9,5}$ \\
     (2,1,2,1,1,2,2,1,3) & $\varphi_{1,5}$ back to (ii)\\
     (2,2,1,2,2,1,2,2,1) & Not needed \\
     
     \hline
     \multicolumn{2}{|c|}{Sections: $D_1$ and $D_4$} \\
     \hline
     Toric pair, cycle of integers & Elementary transformations \\ \hline
     (2,1,1,2,1,2,3,1,2) & $\varphi_{6,2}$, back to (ii) \\
     (2,1,1,2,1,3,1,3,1) & $\varphi_{8,3}$, back to (i)\\
     (2,1,1,1,1,3,2,1,3) & $\varphi_{6,3},\varphi_{9,4}$\\
     (iii) (1,1,1,2,1,3,2,1,3) & $\varphi_{9,3},\varphi_{6,1}$ \\
     (iv) (2,1,1,1,1,4,1,2,2) & $\varphi_{6,3},\varphi_{1,5},\varphi_{6,2}$ \\
     (v)(1,1,1,2,1,4,1,2,2) & $\varphi_{6,2}$, back to (i) \\
     (vi) (1,1,1,1,3,1,3,1,3) & $\varphi_{9,3},\varphi_{5,1}$, back to (i) \\
    
     \hline
     \multicolumn{2}{|c|}{Sections: $D_1$ and $D_3$} \\
     \hline
     Toric pair, cycle of integers & Elementary transformations \\ \hline
     (2,0,3,1,2,2,2,2,1) & $\varphi_{5,9},\varphi_{8,4},\varphi_{3,7}$ \\
     (0,0,1,2,2,2,2,1,5) & $\varphi_{9,2},\varphi_{4,1},\varphi_{9,3},\varphi_{5,1},\varphi_{9,4}$ \\
     (1,0,1,1,3,2,2,1,4) & $\varphi_{5,2}$, back to (iv) \\
     (0,0,1,3,1,3,2,1,4) & $\varphi_{9,2},\varphi_{4,1}$, back to (iii) \\
     (1,0,1,2,1,4,2,1,3) & $\varphi_{6,2},\varphi_{6,2}$, back to (i) \\
     (2,0,1,1,2,3,2,1,3) & $\varphi_{5,2},\varphi_{6,3},\varphi_{9,4}$ \\
     (0,0,1,2,3,1,3,1,4) & $\varphi_{9,2},\varphi_{4,1}$back to (vi)\\
     (1,0,1,1,4,1,3,1,3) & $\varphi_{5,2}$, back to (vi) \\
     (1,0,0,3,2,1,4,1,3) & $\varphi_{4,2},\varphi_{9,3}$, back to (v) \\
     (1,0,1,2,2,1,5,1,2) & $\varphi_{7,2}$, back to (v) \\
     (1,0,2,1,3,1,4,1,2) & $\varphi_{7,2},\varphi_{7,2}$, back to (i) \\
     (1,0,2,2,1,3,3,1,2) & $\varphi_{7,2},\varphi_{5,2}$ \\
     (2,0,2,1,2,2,3,1,2) & $\varphi_{6,2}$, back to (i) \\
     (1,0,0,2,2,3,1,2,4) & $\varphi_{4,2},\varphi_{9,3},\varphi_{6,2}$ \\
     (1,0,1,1,3,3,3,2,3) & $\varphi_{5,2},\varphi_{9,3},\varphi_{5,2}$ \\
     (1,0,2,1,2,4,1,2,2) & $\varphi_{6,2}\varphi_{6,2}$ \\
     (1,0,0,2,3,2,1,3,3) & $\varphi_{4,2},\varphi_{9,3},\varphi_{5,2},\varphi_{8,3}$ \\
     (1,0,1,1,4,2,1,3,2) & $\varphi_{5,2},\varphi_{8,3}$ \\
     (1,0,0,3,2,2,1,2,4) & $\varphi_{4,2},\varphi_{9,3},\varphi_{8,2},\varphi_{4,9},\varphi_{8,3}$ \\
     (1,0,1,2,2,2,1,5,1) & $\varphi_{8,2},\varphi_{4,9},\varphi_{8,3}$ \\
     (1,0,2,1,3,2,1,4,1) & $\varphi_{8,2},\varphi_{5,9}$ \\
     (2,0,2,1,2,3,1,3,1) & $\varphi_{6,2},\varphi_{9,3}$, back to (i) \\

     \hline
    \end{tabular}
    \end{center}\caption{Elementary transformations for $n=9$}\label{table9}
    \end{table}

\end{document}